\documentclass[11pt,a4paper]{article}
\usepackage{a4}
\usepackage{graphicx}
\usepackage{amsmath}
\usepackage{amssymb}
\usepackage[latin1]{inputenc}
\usepackage{here}
\graphicspath{
{./}
{./Figures/}
}

\newcommand{\ub}{{\bar  u}}
\newcommand{\vb}{{\bar  v}}
\newcommand{\wb}{{\bar  w}}
\newcommand{\rhob}{{\bar  \rho}}
\newcommand{\cb}{{\bar  c}}
\newcommand{\euler}{A_{Euler}}
\newcommand{\eulerhathat}{{\hat {\hat A}_{Euler}}}
\newcommand{\eulerhathatpmlone}{{\hat {\hat A}_{Euler}^{pml1}}}
\newcommand{\eulerhathatpmltwo}{{\hat {\hat A}_{Euler}^{pml2}}}
\newcommand{\eulerhathatpmltwoextended}{{\hat {\hat {{\cal A}}}_{Euler}^{pml2}}}
\newcommand{\eulerpmlone}{{ { {{\cal A}}}_{Euler}^{pml1}}}
\newcommand{\eulerpmltwoextended}{{ { {{\cal A}}}_{Euler}^{pml2}}}
\newcommand{\hhp}{\hat {\hat p}}
\newcommand{\hhu}{\hat {\hat u}}
\newcommand{\hhv}{\hat {\hat v}}

\newcommand{\vecstate}{{W}}

\newcommand{\R}{{\bf R}}

\def \ds{\displaystyle}
\def \be{\begin{equation}}
\def \ee{\end{equation}}
\def \bea{\begin{eqnarray}}
\def \eea{\end{eqnarray}}
\def \bean{\begin{eqnarray*}}
\def \eean{\end{eqnarray*}}
\newtheorem{theorem}{Theorem}[section]

\newtheorem{remark}{Remark}

\title{A new construction of perfectly matched layers for the linearized Euler equations}%

\author{
Fr{\'e}d{\'e}ric
Nataf\thanks{CNRS, UMR 7641, CMAP, Ecole Polytechnique, 91128 Palaiseau
Cedex, France}
}

\begin{document}
\maketitle

\begin{abstract} 
Based on a PML for the advective wave equation, we propose two PML models for the linearized Euler equations. The derivation of the first model can be applied to other physical models. The second model was implemented. Numerical results are shown. 
\end{abstract}

\tableofcontents

\bibliographystyle{alpha}

\section{Introduction} 
Since the work by Berenger on perfectly matched layer for the Maxwell equations \cite{MR1294924,MR1412240}  in a computational box, many works have been devoted to a better understanding of their principle and behaviour see \cite{MR1632034}, \cite{Zhao:1996:GTP}, \cite{Chew:1994:3PM}, \cite{Lassas:1998:ECS} \cite{MR1667698} \cite{MR1985305}\cite{MR1916294} \cite{MR1910579} to extensions to other geometries, see \cite{MR2100511} \cite{MR1638033}, or equations see \cite{Hagstrom:2002:ALR} \cite{MR1712596}\cite{MR2077985} . We consider here the linearized Euler equations which has been the subject of many works, see \cite{MR2061242}, \cite{MR1866857}, \cite{MR1633045} \cite{MR1618017} \cite{MR1419743} \cite{MR2077987}  (and references therein). The key issue is a possible lack of long time stability \cite{MR1618017,MR1633045}. A first stable PML was proposed in \cite{MR1866857} for flows normal to the boundary. It has been recently extended to oblique flows in \cite{MR2077987}, see also \cite{Dubois:2002:LTS}. In these works, the PML is obtained via a change of coordinate in the complex plane applied to the direction normal to the boundary (say $x$). This amounts to replacing all the $\partial_x$ derivatives in the Euler system by an operator denoted by $\partial_x^{pml}$ which is still differential in the $x$ direction but pseudo-differential in the other variables. Based on the Smith factorization \cite{Wloka:1995:BVP}, we propose here another strategy for constructing stable PMLs. The main feature is that not all $\partial_x$ derivatives are modified. The Euler system is modified only in order to damp the modes that could reflect at the interface between the physical domain and the PML. In our approach, the vorticity modes which are convected by a transport operator are not damped since they leave naturally the computational domain at an outflow boundary. As a result, even after discretization, the vorticity in the PML region equals (up to the machine accuracy) the vorticity of the reference solution computed in a very large domain, see figure~\ref{fig:omegauntouched} in \S~\ref{sec:numericalresults}.\\
 Moreover, the technique introduced in this paper enables an appropriate ``PML'' treatment of the various scalar equations at the basis of the system of PDEs under consideration. We think therefore that the technique introduced in this paper should extend to other systems of PDEs as well (e.g. shallow water, anistropic elasticity, $\ldots$). \\

More precisely, in section~\ref{sec:smith} we introduce the Smith factorization of the Euler equation. This powerful tool is used in section~\ref{sec:pml} to propose two ways to design PML for the Euler equation (see \cite{Dolean:2005:NCD} for the application of the Smith factorization in domain decomposition methods). Both PMLs are based on the use of a PML for the underlying advective wave equation. The derivation of the first model can be applied to other physical models. In section~\ref{sec:numericalresults}, numerical results are shown for the second model which is easier to implement.

\section{Analysis of the Euler system via Smith factorization}\label{sec:smith}
We write the linearized Euler equations around a constant subsonic flow $(\ub,\vb)$, a constant density $\rhob>0$ and a constant speed of sound $\cb>0$ as:
\be\label{eq:euler}
\left(\begin{array}{ccc} \partial_t + \ub \partial_x + \vb \partial_y &  \rhob \cb^2 \partial_x &   \rhob \cb^2 \partial_y \\
                         \frac{1}{\rhob}\partial_x & \partial_t + \ub \partial_x + \vb \partial_y  & 0 \\
                      \frac{1}{\rhob}\partial_y & 0 & \partial_t +\ub \partial_x+\vb \partial_y 
\end{array} \right)
\ \left( \begin{array}{c} p\\ u\\ v \end{array} \right)
= \ \left( \begin{array}{c} f_p\\ f_u\\ f_v \end{array} \right)
\ee
where $(f_p,f_u,f_v)^T$ is a given right handside.
\subsection{Smith factorization}
We first recall the definition of the Smith factorization of a matrix with polynomial entries and apply it to systems of PDEs, see \cite{Gantmacher:1966:TM1,Gantmacher:1998:TM} (or \cite{Wloka:1995:BVP} for the mathematical analysis of systems of PDEs) and references therein. \\
\begin{theorem}
Let $n$ be an integer and $A$ an invertible $n\times n$ matrix with polynomial entries in one variable $\lambda$ : $A=(a_{ij}(\lambda))_{1\le i,j\le n}$. \\
Then, there exist three matrices with polynomial entries $E$, $D$ and $F$ with the following properties:
\begin{itemize}
\item det($E$) and det($F$) are real numbers.
\item $D$ is a diagonal matrix.
\item $A=EDF$.
\end{itemize}
Morevoer, $D$ is uniquely defined up to a reordering and multiplication of each entry by a constant by  a formula defined as follows. Let $1\le k\le n$,
\begin{itemize}
\item  $S_k$ is the set of all the submatrices of order $k\times k$ extracted from $A$.
\item $Det_k=\{\text{Det}(B_k)\backslash B_k\in S_k  \}$
\item $LD_k$ is a largest common divisor of the set of polynomials $Det_k$. 
\end{itemize}
Then,
\be
D_{kk}(\lambda)=\frac{LD_k(\lambda)}{LD_{k-1}(\lambda)}, \hfill 1\le k\le n
\ee
(by convention, $LD_0=1$). 
\end{theorem}

\paragraph{Application to the Euler system}
We first take formally the Fourier transform of the system (\ref{eq:euler}) with respect to $y$ and $t$ (dual variables are $k$ and $\omega$ resp.). We keep the partial derivatives in $x$ since in the sequel we shall design a PML for a truncation of the domain in the $x$ direction. We note
\begin{equation}\label{eq:eulerhathat}
\eulerhathat = 
\left(\begin{array}{ccc} i\omega + \ub \partial_x + ik \vb  &  \rhob \cb^2 \partial_x &   i \rhob \cb^2 k \\
                         \frac{1}{\rhob}\partial_x & i\omega + \ub \partial_x + ik\vb  & 0 \\
                      \frac{ik}{\rhob} & 0 & i\omega +\ub \partial_x+ i\vb k 
\end{array} \right)
\end{equation}

We can perform a Smith factorization of  $\eulerhathat$ by considering it as a matrix with polynomials in $\partial_x$ entries. We have  
\begin{equation}\label{eq:smithfactorization}
\eulerhathat=EDF
\end{equation}
 where
\begin{equation}\label{eq:Dsmith}
D=\left(\begin{array}{ccc} 1 & 0 & 0 \\
                           0 & 1 & 0 \\
                           0 & 0 & \ds \hat{\hat{\cal G}}\hat{\hat{\cal L}} 
        \end{array}\right) 
\end{equation}
and
\[
E= \frac{1}{(\ub(\cb^2-\ub^2))^{1/3}}\left(\begin{array}{ccc} i \rhob\cb^2k & 0 & 0 \\
                           0 & \ub & 0 \\
       i\omega +\ub\partial_x+i \vb k & E_2 &\ds \frac{\cb^2-\ub^2}{ik\rhob\cb^2}
        \end{array}\right)
\]
and
\[
F= -\left(\begin{array}{ccc} 
\ds \frac{i\omega+\ub\partial_x+ik\vb}{ik\rhob \cb^2} & \ds\frac{\partial_x}{ik} & 1\\
\ds \frac{\partial_x}{\rhob\ub} &\ds \frac{i\omega+\ub\partial_x+ik\vb}{\ub} & 0\\
\ds \frac{\ub}{i\omega+ik\vb} & \ds\frac{\rhob\ub^2}{i\omega+ik\vb} & 0
   \end{array}\right)
\]
where
\[
E_2=\ub\frac{(-\ub\cb^2+\ub^3)\partial_{xx}+(2\ub^2-\cb^2)(i\omega+ik\vb)\partial_x+\ub((i\omega+ik\vb)^2+k^2\cb^2)}{\cb^2(i\omega+ik\vb)},
\]
\begin{equation}
  \label{eq:Grond}
{\hat{\hat {\cal G}}}=i\omega+\ub\partial_x+ik\vb
\end{equation}
 and 
\begin{equation}
  \label{eq:Lrond}
  {\hat{\hat {\cal L}}}= -\omega^2+2ik \ub\vb\partial_{x}+2i\omega(\ub\partial_x+ik\vb)+(\cb^2-\vb^2)k^2  - (\cb^2-\ub^2)\partial_{xx}
\end{equation}
The operators showing up in the diagonal matrix have a physical meaning:
\[
{\cal G}=\partial_t +\ub\partial_x+\vb\partial_y
\] 
is a first order transport operator and 
\[
{\cal L}= \partial_{tt}+2 \ub\vb\partial_{xy}+2\partial_t(\ub\partial_x+\vb\partial_y)-(\cb^2-\vb^2)\partial_{yy}  - (\cb^2-\ub^2)\partial_{xx}
\] 
is the advective wave operator.

\subsection{Modes via Smith factorization}\label{sec:modes}
In the PML analysis of section~\ref{sec:pml}, we shall use the expression of solutions to the homogeneous Euler equation. In order to illustrate the previous section, we make use of the Smith factorization to compute them. We take the Fourier transform in $t$ and $y$ of (\ref{eq:euler}) and seek non zero solutions to 
\[
\eulerhathat \left( \begin{array}{c} \hhp(\omega,x,k)\\ \hhu(\omega,x,k)\\ \hhv(\omega,x,k) \end{array} \right)=0\ \ \ \  x\in\R,\ \omega\in\R,\ k\in\R
\]
Using Smith factorization (\ref{eq:smithfactorization}), we have
\[
EDF \left( \begin{array}{c} \hhp(\omega,x,k)\\ \hhu(\omega,x,k)\\ \hhv(\omega,x,k) \end{array} \right)=0\ \ \ \  x\in\R,\ \omega\in\R,\ k\in\R
\]
Since det$(E)$ is a real number, $E^{-1}$ is still a matrix with polynomials in $\partial_x$ entries so that we can apply it to the above equation and get:
\[
\left(\begin{array}{ccc} 1 & 0 & 0 \\
                           0 & 1 & 0 \\
                           0 & 0 & \ds \hat{\hat{\cal G}}\hat{\hat{\cal L}} 
        \end{array}\right) 
F \left( \begin{array}{c} \hhp(\omega,x,k)\\ \hhu(\omega,x,k)\\ \hhv(\omega,x,k) \end{array} \right)=0\ \ \ \  x\in\R,\ \omega\in\R,\ k\in\R
\]
This implies that 
\begin{equation}\label{eq:eqinF}
F \left( \begin{array}{c} \hhp(\omega,x,k)\\ \hhu(\omega,x,k)\\ \hhv(\omega,x,k) \end{array} \right)
= \left( \begin{array}{c} 0 \\ 0 \\ \sum_i \alpha_i(\omega,k) e^{\lambda_i(\omega,k)x} \end{array} \right)
\ \ \ \  x\in\R,\ \omega\in\R,\ k\in\R
\end{equation}
where $\hat{\hat{{\cal G}}}\hat{\hat{{\cal L}}}(e^{\lambda_i(\omega,k)x})=0$. Since ${\cal G}{\cal L}$ is of third order in the $x$ direction, we have three possible values for $\lambda_i$:
  \begin{equation}
    \label{eq:lambda1}
\ds \lambda_1=-\frac{i\omega+ik\vb}{\ub}    
  \end{equation}
 \begin{equation}
    \label{eq:lambda2}
\ds \lambda_2=
\left\{  \begin{array}{l}       
\ds\frac{u(i\omega+ik\vb)-\cb (i\omega+ik\vb) \sqrt{1-\frac{k^2(\cb^2-\vb^2)}{(\omega+k\vb)^2}}}{\cb^2-\ub^2} \ \ \text{ for }|k|\sqrt{\cb^2-\vb^2}<|\omega+k\vb|        \\  \\
\ds\frac{u(i\omega+ik\vb)-\cb \sqrt{k^2(\cb^2-\vb^2)-(\omega+k\vb)^2)}}{\cb^2-\ub^2} \ \ \text{ for }|k|\sqrt{\cb^2-\vb^2}>|\omega+k\vb|
\end{array}\right.
  \end{equation}

 \begin{equation}
    \label{eq:lambda3}
\ds \lambda_3=
\left\{  \begin{array}{l}       
\ds\frac{u(i\omega+ik\vb)+\cb (i\omega+ik\vb) \sqrt{1-\frac{k^2(\cb^2-\vb^2)}{(\omega+k\vb)^2}}}{\cb^2-\ub^2} \ \ \text{ for }|k|\sqrt{\cb^2-\vb^2}<|\omega+k\vb|        \\  \\
\ds\frac{u(i\omega+ik\vb)+\cb \sqrt{k^2(\cb^2-\vb^2)-(\omega+k\vb)^2)}}{\cb^2-\ub^2} \ \ \text{ for }|k|\sqrt{\cb^2-\vb^2}>|\omega+k\vb|
\end{array}\right.
  \end{equation}
\begin{remark}
$\lambda_1$ comes from the transport operator ${\cal G}$ whereas $\lambda_{2,3}$ come from the advective wave operator ${\cal L}$. 
\end{remark}
Since det$(F)$ is a real number, $F^{-1}$  is still a matrix with polynomials in $\partial_x$ entries so that we can apply it to equation (\ref{eq:eqinF}) and get:
\begin{equation}
   \left( \begin{array}{c} \hhp(\omega,x,k)\\ \hhu(\omega,x,k)\\ \hhv(\omega,x,k) \end{array} \right)
= \sum_{i=1}^3 \alpha_i(\omega,k) F^{-1}
 \left( \begin{array}{c} 0 \\ 0 \\ e^{\lambda_i(\omega,k)x} \end{array} \right)
\end{equation}
We shall define, for $i=1,2,3$ 
\begin{equation}\label{eq:modeseuler}
W_i(\omega,k)= e^{-\lambda_i(\omega,k)x}  F^{-1}
 \left( \begin{array}{c} 0 \\ 0 \\ e^{\lambda_i(\omega,k)x} \end{array} \right)
\end{equation}
It is easy to check that indeed $W_i$ does not depend on $x$ and that
\begin{equation}
  \label{eq:modes}
   \left( \begin{array}{c} \hhp(\omega,x,k)\\ \hhu(\omega,x,k)\\ \hhv(\omega,x,k) \end{array} \right)
= \sum_{i=1}^3 \alpha_i(\omega,k) W_i(\omega,k) e^{\lambda_i(\omega,k)x}
\end{equation}
In the classical mode analysis of a system of PDEs, vectors $W_i(\omega,k)$ are obtained after the diagonalization of a matrix. Using the Smith factorization simplifies the computation since these vectors are given by the explicit formula \eqref{eq:modeseuler} and not by an eigenvalue computation.\\

\noindent{\bf Notation} Let $i$ and $j$ be integers. For a matrix $A$, $(A)_{ij}$ denotes the entry of the $i$-th row and of the $j$-th column. For a vector $V$, $(V)_i$ denotes its $i$-th component.
\begin{remark}
 In section~\ref{sec:pmlness2}, we shall use that 
  \begin{equation}
    \label{eq:Fm113}
    (F^{-1})_{13}=\frac{\hat{\hat{{{\cal G}}}}}{\ub}
  \end{equation}
\end{remark}

\section{PMLs for the Euler System }\label{sec:pml}
The Smith factorization of the Euler system (\ref{eq:smithfactorization}) and the computations of the previous section show that the modes correspond either to operator ${\cal L}$ or to operator ${\cal G}$. Among these two operators, the only operator which generates waves propagating in both positive $x$ and negative $x$ directions is the operator ${\cal L}$. This  suggests that designing a PML  for the Euler equation can be reduced to the design of  PML for the advective wave operator ${\cal L}$.

\subsection{PML for the advective wave equation}
\label{sec:pml-wave-eq}
This question has been the subject of several works, \cite{MR2077987}, \cite{Dubois:2002:LTS},  \cite{Hagstrom:2002:ALR},   \cite{MR2051072} \cite{MR2077985} \cite{MR2077990} and references therein. Following these works, we use the for operator ${\cal L}$ a PML defined by replacing the $x$ derivatives by a ``pml'' $x$ derivative. The definition is as follows. Let $\sigma(\omega,x,k)\ge 0$ be the damping parameter of the PML and ${\cal F}^{-1}$ is the inverse Fourier transform in the variables $\omega$ and $k$, we define the pseudo-differential operator  $\alpha(x)$ :
\begin{equation}
  \label{eq:alphax}
  \alpha(x)(\phi)={\cal F}^{-1}(\frac{\cb(i\omega+ik\vb)}{\cb(i\omega+ik\vb)+(\cb^2-\ub^2)\sigma(\omega,x,k)}\,\,\hat{\hat{\phi}})
\end{equation}
We define the $\partial_x^{pml}$ derivative by 
\begin{equation}
  \label{eq:dxpml}
  \partial_x^{pml}=\alpha(x)[\partial_x-\frac{\ub}{\cb^2-\ub^2}(\partial_t+\vb\partial_y)]+\frac{\ub}{\cb^2-\ub^2}(\partial_t+\vb\partial_y)
\end{equation}
We are now in position to define the PML-$x$ equations of the advective wave operator $\mathcal{L}$:
\begin{equation}
  \label{eq:Lrondpml}
  {\cal L}_{pml}= \partial_{tt}+2 \ub\vb\partial_y(\partial_{x}^{pml})+2\partial_t(\ub\partial_x^{pml}+\vb\partial_y)-(\cb^2-\vb^2)\partial_{yy}  - (\cb^2-\ub^2)(\partial_{x}^{pml})^2
\end{equation}
Let us notice that we have
\begin{equation}
  \label{eq:dxpmlmdx}
  \partial_x^{pml}-\partial_x=\gamma(x)[\partial_x-\frac{\ub}{\cb^2-\ub^2}(\partial_t+\vb\partial_y)]
\end{equation}
where  the operator $\gamma(x)$ is a pseudo-differential operator in the $t$ and $y$ variables: 
\begin{equation}
  \label{eq:gammaaux}
    \gamma(x)(\phi)={\cal F}^{-1}(\frac{-(\cb^2-\ub^2)\sigma(\omega,x,k)}{\cb(i\omega+ik\vb)+(\cb^2-\ub^2)\sigma(\omega,x,k)}\,\,\hat{\hat{\phi}})
\end{equation}
A PML-$y$ used for truncating the domain in the $y$ direction would consist in replacing in the operator ${\cal L}$ the $y$ derivatives by a ``pml'' derivative in the $y$ direction defined as follows:
\begin{equation}
  \label{eq:dypml}
  \partial_y^{pml}=\alpha(y)[\partial_y-\frac{\vb}{\cb^2-\vb^2}(\partial_t+\ub\partial_x)]+\frac{\vb}{\cb^2-\vb^2}(\partial_t+\ub\partial_x)
\end{equation}
The PML-$y$ equations of the advective wave operator $\mathcal{L}$ for the truncation of the domain in the $y$ direction read:
\begin{equation}
  \label{eq:Lrondpmly}
  {\cal L}_{pml,y}= \partial_{tt}+2 \ub\vb\partial_y^{pml}(\partial_{x})+2\partial_t(\ub\partial_x+\vb\partial_y^{pml})-(\cb^2-\vb^2)(\partial_{y}^{pml})^2  - (\cb^2-\ub^2)\partial_{xx}
\end{equation}
In order to give a complete definition of a PML bordering a rectangular computational domain, we have three possibilities for the corner region.The first one consists in desiging a third PML model in the corner that is compatible with both PML-$x$ and PML-$y$ as was done for the Maxwell system in \cite{MR1294924} for instance. The second possibility is to use prismatoidal coordinates \cite{Lassas:2001:APE}. The advantage is that it allows for arbitrary convex computational domains and not only rectangular ones. The third possibility consists simply in placing side by side PML-$x$ and PML-$y$ regions. We have implemented this last simple approach and obtain good results, see \S~\ref{sec:numericalresults}. Of course, the two first options deserve further investigations. 

\subsection{First PML model}\label{sec:firstpml}
Based on (\ref{eq:smithfactorization}), a first possibility is to define a PML for the Euler system by substitution of  ${\cal L}$ with ${\cal L}^{pml}$ in matrix $D$ (see formula (\ref{eq:Dsmith})).  In matrices $E$ and $F$ and in the operator ${\cal G}$, the $x$ derivatives are not modified. We modify only the advective wave operator. Let
\begin{equation}
  \label{eq:Dpml1}
   D^{pml}= \left(\begin{array}{ccc} 1 & 0 & 0 \\
                           0 & 1 & 0 \\
                           0 & 0 & \ds \hat{\hat{{\cal G}}}\hat{\hat{{\cal L}}}^{pml}
        \end{array}\right)
\end{equation}
we define a PML for the Euler system
\begin{equation}
  \label{eq:eulerpml1}
  \eulerhathatpmlone=E D^{pml}  F
\end{equation}
with the following interface conditions between the Euler media and the PML
\begin{equation}
  \label{eq:ICforPML1}
  \begin{array}{c}
  {\mathcal L}((F(\vecstate_{Euler}))_3)= {\mathcal L}_{pml}((F(\vecstate_{pml}))_3)\\
  {\mathcal G}((F(\vecstate_{Euler}))_3)= {\mathcal G}((F(\vecstate_{pml}))_3)\\
  \partial_x({\mathcal G}((F(\vecstate_{Euler}))_3))= \partial_x^{pml}({\mathcal G}((F(\vecstate_{pml}))_3))
  \end{array}
\end{equation}
where the subscript ${}_3$ denotes the third component of a vector.
\paragraph{Study of the first PML media}
We now give the results of an analysis of the PML system similar to that of \S~\ref{sec:modes} for the Euler system. We define
\begin{equation*}
\left\{  \begin{array}{l}
   \lambda_1^{pml}=\lambda_1\\
  \ds  \lambda_{2,3}^{pml}= 
(1+\frac{(\cb^2-\ub^2)\sigma}{\cb(i\omega+ik\vb)})
(\lambda_{2,3}-\frac{\ub}{\cb^2-\ub^2}(i\omega+i\vb k))
+\frac{\ub}{\cb^2-\ub^2}(i\omega+i\vb k)
  \end{array}\right.
\end{equation*}
and for $i=1,2,3$ 
\begin{equation}\label{eq:vectorspml1}
W_i^{pml1}(\omega,k)= e^{-\lambda_i^{pml}(\omega,k)x}  F^{-1}
 \left( \begin{array}{c} 0 \\ 0 \\ e^{\lambda_i^{pml}(\omega,k)x} \end{array} \right)
\end{equation}
It is easy to check that indeed $W_i^{pml1}$ does not depend on $x$ and that for any solution $W^{pml1}$ to the homogeneous first PML system there exist $(\beta_i(\omega,k))_{1\le i\le 3}$ such that
\begin{equation}
  \label{eq:modespml1}
{\hat {\hat W}}^{pml1}(\omega,x,k)
= \sum_{i=1}^3 \beta_i(\omega,k) W_i^{pml1}(\omega,k) e^{\lambda_i^{pml1}(\omega,k)x}
\end{equation}
If we consider the solution in the positive $x$ half space, its  boundedness as $x$ tends to infinity implies that $\beta_3=0$.\\

\subsection{PMLness of the first model}\label{sec:pmlness1}
A key property of a PML is that there is no reflection at the interface between the Euler media and the PML media. We will prove that it is the case for a truncation of the space with an infinite PML starting at $x=0$ with a $x$ independant damping parameter $\sigma$. We have to consider the following coupled problem:\\
 Find $(W_{l},W_{r})=((p_l,u_l,v_l),(p_r,u_r,v_r))$ such that:
\begin{equation}\label{eq:coupledsysteme1}
  \begin{array}{l}
\euler W_l =  0, \ \ t>0,\ x<0,\ y\in\R\\  \\
 \eulerpmlone W_r = 0,  \ \ t>0,\ x>0,\ y\in\R 
\end{array}
\end{equation}
with the interface conditions~\eqref{eq:ICforPML1}
\begin{equation}\label{eq:ICforPML1bis}
  \begin{array}{c}
  {\mathcal L}((F(W_l))_3)= {\mathcal L}_{pml}((F(W_r))_3)\\
  {\mathcal G}((F(W_l))_3)= {\mathcal G}((F(W_r))_3)\\
  \partial_x({\mathcal G}((F(W_l))_3))= \partial_x^{pml}({\mathcal G}((F(W_r))_3))
  \end{array}
\end{equation}
We take the Fourier transform in $t$ and $y$ of \eqref{eq:coupledsysteme1} and get:
\begin{equation}\label{eq:coupledsystemeFourier1}
  \begin{array}{l}
\eulerhathat \hat{\hat{W_l}} =  0, \ \  x<0,\ \omega,\,k\in\R\\  \\
\eulerhathatpmlone  \hat{\hat{W_r}} = 0,  \  x>0,\ \omega,\,k\in\R 
\end{array}
\end{equation}
From section~\ref{sec:modes}, we know that the general solution to the Euler system is :
\begin{equation}\label{eq:WLpmlness1}
\hat{\hat{W}}_l=\sum_{i=1}^3 \alpha_i(\omega,k) W_i(\omega,k) e^{\lambda_i(\omega,k)x}
\end{equation}
where $W_i$ is defined in (\ref{eq:modeseuler}). As for the solution in the PML media, we know from~\eqref{eq:modespml1} and the boundedness of the solution as $x$ tends to infinity that 
\begin{equation}\label{eq:WRpmlness1}
\hat{\hat{W}}_r=\sum_{i=1}^2 \beta_i(\omega,k) W_i^{pml1}(\omega,k) e^{\lambda_i^{pml}(\omega,k)x}
\end{equation}
We study the adequacy of the PML by considering $\alpha_1$ and $\alpha_2$ to be given. This corresponds to ingoing waves from the Euler media and moving towards the interface between the Euler media and the PML media. The three other quantities ($\alpha_3,(\beta_i)_{i=1,\ldots,2})$ are determined by the interface conditions~\eqref{eq:ICforPML1bis}. The media is perfectly matched if we have no reflection in the Euler media, i.e. if $\alpha_3=0$.  We now prove that this is indeed the case.\\
From \eqref{eq:WLpmlness1} and \eqref{eq:modeseuler}, we have:
\begin{equation}
  \label{eq:F3WL}
    (F({\hat{\hat W}}_l))_3=\sum_{i=1}^3 \alpha_i(\omega,k) e^{\lambda_i(\omega,k)x}
\end{equation}
From \eqref{eq:WRpmlness1} and \eqref{eq:vectorspml1}, we have:
\begin{equation}
  \label{eq:F3WR}
    (F({\hat{\hat W}}_r))_3=\sum_{i=1}^2 \beta_i(\omega,k) e^{\lambda_i^{pml}(\omega,k)x}
\end{equation}
Let 
\[
a=\frac{\cb(i\omega+ik\vb)}{\cb(i\omega+ik\vb)+(\cb^2-\ub^2)\sigma(\omega,x,k)}
\]
Notice that
\[
\begin{array}{l}
{\hat{\hat {\mathcal G}}}=\ub(\partial_x-\lambda_1)\\
{\hat{\hat {\mathcal L}}}=-(\cb^2-\ub^2)(\partial_x-\lambda_2)(\partial_x-\lambda_3)\\
{\hat{\hat {\mathcal L}}}^{pml}=-(\cb^2-\ub^2)a^2(\partial_x-\lambda_2^{pml})(\partial_x-\lambda_3^{pml})
\end{array}
\]
and that
\[
\partial_x^{pml}(e^{\lambda_i^{pml}x})=\lambda_i\,e^{\lambda_i^{pml}x}
\]
so that the interface conditions \eqref{eq:ICforPML1bis} yield:
\begin{equation}\label{eq:FourierICforPML1bis}
  \begin{array}{c}
\alpha_1 (\lambda_1-\lambda_2)(\lambda_1-\lambda_3)=\beta_1 a^2(\lambda_1^{pml1}-\lambda_2^{pml})(\lambda_1^{pml}-\lambda_3^{pml}) \\ \\
\alpha_2 (\lambda_2-\lambda_1)+\alpha_3  (\lambda_3-\lambda_1)=\beta_2(\lambda_2^{pml}-\lambda_1) \\ \\
\alpha_2 (\lambda_2-\lambda_1)\lambda_2+\alpha_3  (\lambda_3-\lambda_1)\lambda_3=\beta_2(\lambda_2^{pml}-\lambda_1)\lambda_2 
  \end{array}
\end{equation}
Thus we have
\[
\begin{array}{l}
\ds \beta_1=\alpha_1 a^2 \frac{(\lambda_1-\lambda_2)(\lambda_1-\lambda_3)}{(\lambda_1^{pml1}-\lambda_2^{pml})(\lambda_1^{pml}-\lambda_3^{pml})}\\ \\
\ds \beta_2= \alpha_2 \frac{\lambda_2-\lambda_1}{\lambda_2^{pml}-\lambda_1} \\ \\
\alpha_3=0
\end{array}
\]
The nullity of $\alpha_3$ shows that there is no reflection at the interface between the Euler and an infinite PML media. \\
In practice, the PML has a finite width. As a result, the coefficient $\beta_3$ will not be zero in general. But, the corresponding mode $W_3^{pml}\,e^{\lambda_3^{pml}x}$ is exponentially decreasing as $x$ is increasing. Its contribution to the solution on the interface can be made as small as necessary simply by increasing the width of the PML.

\paragraph{Complexity of the first PML model}

A direct computation yields:
\begin{equation}
  \label{eq:diffeulerpml}
    \eulerhathatpmlone=  \eulerhathat+
 \left(\begin{array}{ccc}  0 & 0 & 0 \\
                           0 & 0 & 0 \\
                           C_1 & C_2 & 0
        \end{array}\right)
\end{equation}
where
$$
C_1=\frac{(\partial_x-\partial_x^{pml}) \hat{\hat{{\cal G}}}[(\ub^2-\cb^2)(\partial_x+\partial_x^{pml})+2\ub (i\omega+i\vb k)]}{i\rhob\cb^2 k(i\omega+ik\vb)}
\text{ and }C_2=\frac{C_1}{\rhob \ub}
$$
The difference with the Euler system concerns only the last equation on the variable $v$, but : 
\begin{enumerate}
\item The formula is complex and involves  third order derivatives on both the pressure $p$ and the normal velocity $u$. 
\item The formula implies a division by $i\rhob\cb^2 k(i\omega+ik\vb)$ which can be zero. 
\end{enumerate}
As for the first point, one might argue that it is just a matter of implementation. The second point seems more serious. A possible cure could be to take:
\[
\sigma(\omega,x,k)=\tilde\sigma(x)\left(\rhob\cb^2 k(\omega+k\vb)\right)^2
\]
where $\tilde\sigma(x)\ge 0$. From formulas for $C_1$ and $C_2$ and formula (\ref{eq:dxpmlmdx})-(\ref{eq:gammaaux}), we get then
\[
C_1=\tilde\sigma(x)\frac{(i\rhob\cb^2 k(i\omega+ik\vb))(\cb^2-\ub^2)}{\cb(i\omega+ik\vb)+(\cb^2-\ub^2)\sigma(\omega,x,k)}
 \hat{\hat{{\cal G}}}[(\ub^2-\cb^2)(\partial_x+\partial_x^{pml})+2\ub (i\omega+i\vb k)]
\]

As a result the symbols of $C_1$ and $C_2$ would not vanish. But it would be at the expense of the damping of the PML. Indeed, $\sigma(\omega,x,k)$ would be small for small values of $k$ or of $i\omega+ik\vb$. The present first model raises difficulties. Nevertheless, it should deserve interest since it corresponds to a systematic way to design a PML for systems of PDEs. Moreover, since matrices $E$ and $F$ are not unique, it is quite possible that a more suitable Smith factorization when used in formula~(\ref{eq:eulerpml1}) would lead to a practicable PML. In the next section, we design another PML for the Euler system whose numerical results will be given in section~\ref{sec:numericalresults}. 

\subsection{Second PML model}
The rationale for this model is that the pressure $p$ satisfies an advective wave equation which is the only equation that demands a PML. Indeed, let multiply (\ref{eq:eulerhathat}) by the matrix
\begin{equation}
  \label{eq:Emodel2}
  El= \left(\begin{array}{ccc}  \hat{\hat{\cal G}} & -\rhob\cb^2\partial_x & -i\rhob\cb^2 k \\
                           0 & 1 & 0 \\
                           0 & 0 & 1
        \end{array}\right)
\end{equation}
We get:
\begin{equation}
  \label{eq:eulerpressure}
El\,\eulerhathat=\left(\begin{array}{ccc} \hat{\hat{\cal L}}  & 0 & 0 \\
                         \frac{1}{\rhob}\partial_x & i\omega + \ub \partial_x + ik\vb  & 0 \\
                      \frac{ik}{\rhob} & 0 & i\omega +\ub \partial_x+ i\vb k 
\end{array} \right)
\end{equation}
We substitute $\hat{\hat{\cal L}}$ with ${\hat{\hat{\cal L}}}^{pml}$ and apply 
$$
El^{-1}=\left(\begin{array}{ccc}  \hat{\hat{{\cal G}}}^{-1} & -\rhob\cb^2\partial_x \hat{\hat{{\cal G}}}^{-1}  & -i\rhob\cb^2 k \hat{\hat{{\cal G}}}^{-1}  \\
                           0 & 1 & 0 \\
                           0 & 0 & 1
        \end{array}\right)
$$ and we are thus led to define:
\begin{equation}
  \label{eq:eulerpml2}
  \eulerhathatpmltwo=\left(\begin{array}{ccc} \hat{\hat{{\cal G}}}^{-1} (\hat{\hat{{\cal L}}}^{pml}+\cb^2(\partial_{xx}-k^2))  & \rhob \cb^2 \partial_x & i \rhob \cb^2 k \\
                         \frac{1}{\rhob}\partial_x &  \hat{\hat{{\cal G}}}  & 0 \\
                      \frac{ik}{\rhob} & 0 &  \hat{\hat{{\cal G}}}
\end{array} \right)
\end{equation}
A direct computation yields:
\begin{equation}
  \label{eq:diffpmltwo}
      \eulerhathatpmltwo=  \eulerhathat+
 \left(\begin{array}{ccc}  (\hat{\hat{{\cal L}}}^{pml}-\hat{\hat{{\cal L}}}) \hat{\hat{{\cal G}}}^{-1} & 0 & 0 \\
                           0 & 0 & 0 \\
                           0 & 0 & 0
        \end{array}\right)
\end{equation}
In order to get rid of the operator $ \hat{\hat{{\cal G}}}^{-1}$, we introduce a new variable ${\cal P}$ such that ${\cal G}({\cal P})=p$ so that in the physical space the enlarged PML system we consider reads:
\begin{equation}
  \label{eq:eulerpml2extended}
  \eulerpmltwoextended
 \left( \begin{array}{c}{{\cal P}} \\ p\\ u\\ v \end{array} \right)
=\left(\begin{array}{cccc} 
               {{{\cal G}}}  & -1 & 0 & 0 \\
 {{{\cal L}}}^{pml}-{{{\cal L}}}& {{{\cal G}}} &  \rhob\cb^2 \partial_x  &  \rhob\cb^2 \partial_y\\
              0 &        \frac{1}{\rhob}\partial_x &   {{{\cal G}}} & 0  \\
              0 &    \frac{1}{\rhob}\partial_y & 0 & {{{\cal G}}}  
\end{array} \right)
 \left( \begin{array}{c}{{\cal P}}(t,x,y) \\ p(t,x,y)\\ u(t,x,y)\\ v(t,x,y) \end{array} \right)
=0, \ \ t>0, x>0, y\in \R
\end{equation}
with the following interface conditions between the Euler media and the PML
\[
 {{\cal P}}=0,\ p\text{ and }u \text{ are continuous,}\ \partial_x(p_{Euler})=\partial_x^{pml}(p_{pml})
\]
\paragraph{Study of the second PML media}
We now proceed to an analysis of the PML system similar to that of \S~\ref{sec:modes} for the Euler system. The Smith factorization of $\eulerhathatpmltwoextended$ reads
\[
\eulerhathatpmltwoextended = \tilde E \tilde D \tilde F
\]
where 
\begin{equation}
  \label{eq:diagpmlsmith}
  \tilde D= \left(\begin{array}{cccc} 
               1 & 0 & 0 & 0 \\
               0 & 1 & 0 & 0 \\
              0 & 0 &   \hat{\hat{{\cal G}}} & 0  \\
              0 & 0 & 0 & \hat{\hat{{\cal G}}} \hat{\hat{{\cal L}}}^{pml} 
\end{array} \right)
\end{equation}
and $\tilde E$ and $\tilde F$ are matrices with polynomial in $\partial_x$ entries and their determinants are one. This will enable us to give the general form of the solutions to the homogeneous PML equations. Indeed, let us denote $\hat{\hat{W}}=(\hat{\hat{{{\cal P}}}},\hat{\hat{p}},\hat{\hat{u}},\hat{\hat{v}})^T$ such a solution. From the Smith factorization, there exist $(\beta_i(\omega,k))_{i=0,\ldots,3}$ such that 
\[
\tilde F(\hat{\hat{W}})= \beta_0  \left(\begin{array}{c} 
               0 \\ 0 \\ e^{\lambda_0^{pml} x} \\ 0
\end{array} \right)
+
\sum_{i=1}^3 \beta_i  \left(\begin{array}{c} 
               0 \\ 0 \\ 0 \\ e^{\lambda_i^{pml} x} 
\end{array} \right)
\]
where 
\begin{equation*}
\left\{  \begin{array}{l}
   \lambda_0^{pml}=\lambda_1^{pml}=\lambda_1\\
  \ds  \lambda_{2,3}^{pml}= 
(1+\frac{(\cb^2-\ub^2)\sigma}{\cb(i\omega+ik\vb)})
(\lambda_{2,3}-\frac{\ub}{\cb^2-\ub^2}(i\omega+i\vb k))
+\frac{\ub}{\cb^2-\ub^2}(i\omega+i\vb k)
  \end{array}\right.
\end{equation*}
By applying ${\tilde F}^{-1}$ to the above equation, we see that there exist vectors $W_i^{pml2}(\omega,k)$, $i=0,\ldots,3$ such that 
\begin{equation}\label{eq:wpml2}
\hat{\hat{W}}=\sum_{i=0}^3 \beta_i(\omega,k) W_i^{pml2}(\omega,k) e^{\lambda_i^{pml}(\omega,k)x}
\end{equation}
If we consider the solution in the positive $x$ half space, its  boundedness as $x$ tends to infinity implies that $\beta_3=0$.\\

\subsection{PMLness of the second model}\label{sec:pmlness2}
We proceed similarly to \S~\ref{sec:pmlness1}. We will prove there is no reflection at the interface between the Euler media and the PML media for a truncation of the space with an infinite PML starting at $x=0$ with a constant damping parameter $\sigma$. We consider the following coupled problem:\\
 Find $(W_{l},W_{r})=((p_l,u_l,v_l),({{\cal P}}_r,p_r,u_r,v_r))$ such that:
\begin{equation}\label{eq:coupledsysteme}
  \begin{array}{l}
\euler W_l =  0, \ \ t>0,\ x<0,\ y\in\R\\  \\
 \eulerpmltwoextended W_r = 0,  \ \ t>0,\ x>0,\ y\in\R \\  \\
\text{ at }x=0,\ \ {{\cal P}}_r=0,\ p_l=p_r,\ \partial_x(p_l)=\partial_x^{pml}(p_r),\ u_l=u_r, \ \ t>0,\ y\in\R
\end{array}
\end{equation}
We take the Fourier transform in $t$ and $y$ of the above coupled system and get:
\begin{equation}\label{eq:coupledsystemeFourier}
  \begin{array}{l}
\eulerhathat \hat{\hat{W_l}} =  0, \ \  x<0,\ \omega,\,k\in\R\\  \\
\eulerhathatpmltwoextended  \hat{\hat{W_r}} = 0,  \  x>0,\ \omega,\,k\in\R \\  \\
\text{ at }x=0,\ \ \hat{\hat{{{\cal P}}}}_r=0,\ \hat{\hat{p}}_l=\hat{\hat{p}}_r,\ \partial_x(\hat{\hat{p}}_l)=\partial_x^{pml}(\hat{\hat{p}}_r),\ \hat{\hat{u}}_l=\hat{\hat{u}}_r, \ \omega,\,k\in\R
\end{array}
\end{equation}

From section~\ref{sec:modes}, we know that the general solution to the Euler system is :
\[
\hat{\hat{W}}_l=\sum_{i=1}^3 \alpha_i(\omega,k) W_i(\omega,k) e^{\lambda_i(\omega,k)x}
\]
where $W_i$ is defined in (\ref{eq:modeseuler}). As for the solution in the PML media, we know from (\ref{eq:wpml2}) and the boundedness of the solution as $x$ tends to infinity that 
\[
\hat{\hat{W}}_r=\sum_{i=0}^2 \beta_i(\omega,k) W_i^{pml2}(\omega,k) e^{\lambda_i^{pml}(\omega,k)x}
\]
We study the adequacy of the PML by considering $\alpha_1$ and $\alpha_2$ to be given. This corresponds to ingoing waves from the Euler media and moving towards the interface between the Euler media and the PML media. The four other quantities ($\alpha_3,(\beta_i)_{i=0,\ldots,2})$ are determined by the interface conditions. The media is perfectly matched if we have no reflection in the Euler media, i.e. if $\alpha_3=0$.  We now prove that this is indeed the case. We focuse on the equation satisfied by the pressure. By applying matrix 
\[
\left(\begin{array}{cccc}  \hat{\hat{\cal G}} & \hat{\hat{\cal G}} & -\rhob\cb^2\partial_x & -i\rhob\cb^2 k 
        \end{array}\right)
\]
to $\eulerhathatpmltwoextended$, we have that $\hat{\hat{p}}_r$ satisfies the equation of advective Helmholtz PML media:
\[
         \hat{\hat{{{\cal L}}}}^{pml}(\hat{\hat{p}}_r)=0
\]
From (\ref{eq:eulerpressure}), we also have that 
\[
         \hat{\hat{{{\cal L}}}}(\hat{\hat{p}}_l)=0
\]
The interface conditions on the  pressure are 
$\hat{\hat{p}}_l=\hat{\hat{p}}_r$ and  $\partial_x(\hat{\hat{p}}_l)=\partial_x^{pml}(\hat{\hat{p}}_r)$. We know from works on PML for the convective Helmholtz that there is no reflection at the interface for the pressure. Therefore there exists $\beta_p(\omega,k)$ such that 

\begin{equation}\label{eq:pl}
\hat{\hat{p}}_l = \beta_p(\omega,k)e^{\lambda_2(\omega,k)x}
\end{equation}
We prove now that as a consequence, $\alpha_3$ is zero. Indeed, taking the first component of (\ref{eq:wpml2}) we have that 
\begin{equation}\label{eq:plviaw}
\hat{\hat{p}}_l = \sum_{i=1}^3 \alpha_i\, (W_i)_1\, e^{\lambda_i x}
\end{equation}
From (\ref{eq:modeseuler}) and (\ref{eq:Fm113}), we have
\[
(W_i)_1=\frac{1}{\ub}\hat{\hat{{\cal G}}}(e^{\lambda_i x})=\frac{1}{\ub} (\lambda_i-\lambda_1)\,e^{\lambda_i x}
\]
So that we have $(W_1)_1=0$ whereas $(W_i)_1\neq 0$ for $i=2,\,3$. Thus we can infer from (\ref{eq:pl}) and (\ref{eq:plviaw}) that $\alpha_3=0$. This shows that there is no reflection at the interface between the Euler and an infinite PML media. \\
In practice, the PML has a finite width. As a result, the coefficient $\beta_3$ will not be zero in general. But, the corresponding mode $W_3^{pml}\,e^{\lambda_3^{pml}x}$ is exponentially decreasing as $x$ is increasing. Its contribution to the solution on the interface can be made as small as necessary simply by increasing the width of the PML.

\begin{remark}[PMLs for the 3D Euler system]
Using obvious notations, we briefly mention that the design of a PML for the three-dimensional Euler system:
\be\label{eq:euler3D}
\left(\begin{array}{cccc} \partial_t + \ub \partial_x + \vb \partial_y + \wb \partial_z &  \rhob \cb^2 \partial_x &   \rhob \cb^2 \partial_y  &   \rhob \cb^2 \partial_z \\
                         \frac{1}{\rhob}\partial_x & \partial_t + \ub \partial_x + \vb \partial_y + \wb \partial_z & 0 & 0 \\
                      \frac{1}{\rhob}\partial_y & 0 & \partial_t +\ub \partial_x+\vb \partial_y + \wb \partial_z & 0 \\
                      \frac{1}{\rhob}\partial_z & 0 & 0 & \partial_t +\ub \partial_x+\vb \partial_y + \wb \partial_z
\end{array} \right)
\ee
is very similar to the two-dimensional case. Let $\mathcal{G}_3$ be the 3D first order transport operator and $\mathcal{L}_3$ the 3D acoustic wave equation. The Smith form of the 3D compressible Euler equations is
\begin{equation*}
     D=
 \left(\begin{array}{cccc} 1 & 0 & 0 & 0 \\
                           0 & 1 & 0 & 0\\
                           0 & 0 & \mathcal{G}_3 & 0\\
                           0 & 0 & 0 &  \mathcal{G}_3\mathcal{L}_3
        \end{array}\right)
\end{equation*}
As in the 2D case, only the advective wave operator $\mathcal{L}_3$ needs a ``pml'' procedure. For instance, the second model should read:
\begin{equation*}
\left(\begin{array}{crccc} 
               {{{\cal G}_3}}  & -1 & 0 & 0 & 0\\
 {{{\cal L}}}^{pml}_3-{{{\cal L}}}_3& {{{\cal G}}}_3 &  \rhob\cb^2 \partial_x  &  \rhob\cb^2 \partial_y &  \rhob\cb^2 \partial_z\\
              0 &      \frac{1}{\rhob}\partial_x &   {{{\cal G}}} & 0 & 0 \\
              0 &    \frac{1}{\rhob}\partial_y & 0 & {{{\cal G}}} & 0 \\  
              0 &    \frac{1}{\rhob}\partial_z & 0 & 0 & {{{\cal G}}}  
\end{array} \right)
 \left( \begin{array}{c}{{\cal P}} \\ p\\ u\\ v \\ w \end{array} \right)
=0
\end{equation*}
\end{remark}

\section{Numerical Results}\label{sec:numericalresults}
\begin{figure}
  \centering
  \includegraphics[width=0.48\textwidth,height=0.35\textwidth]
    {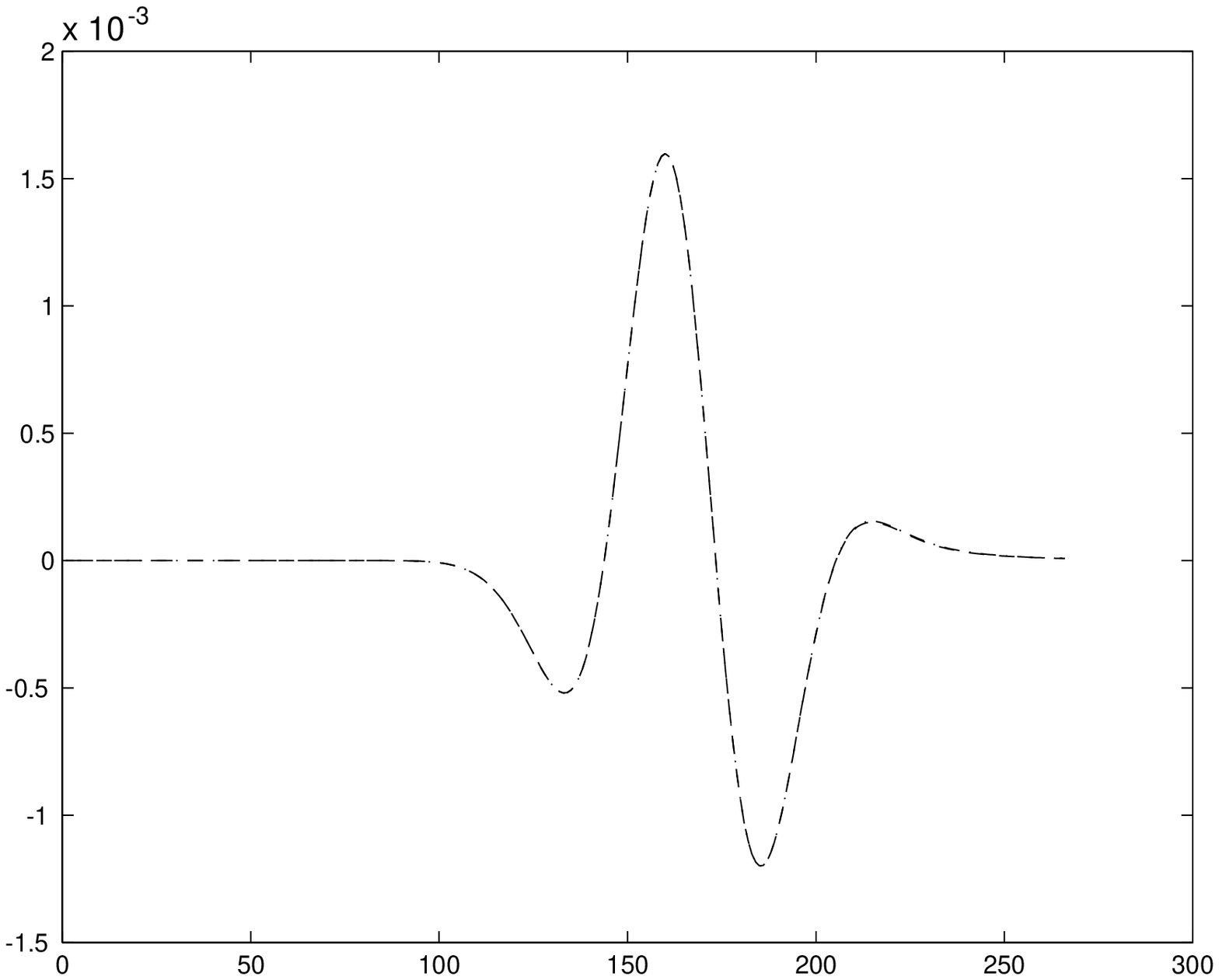}
  \includegraphics[width=0.48\textwidth,height=0.35\textwidth]
    {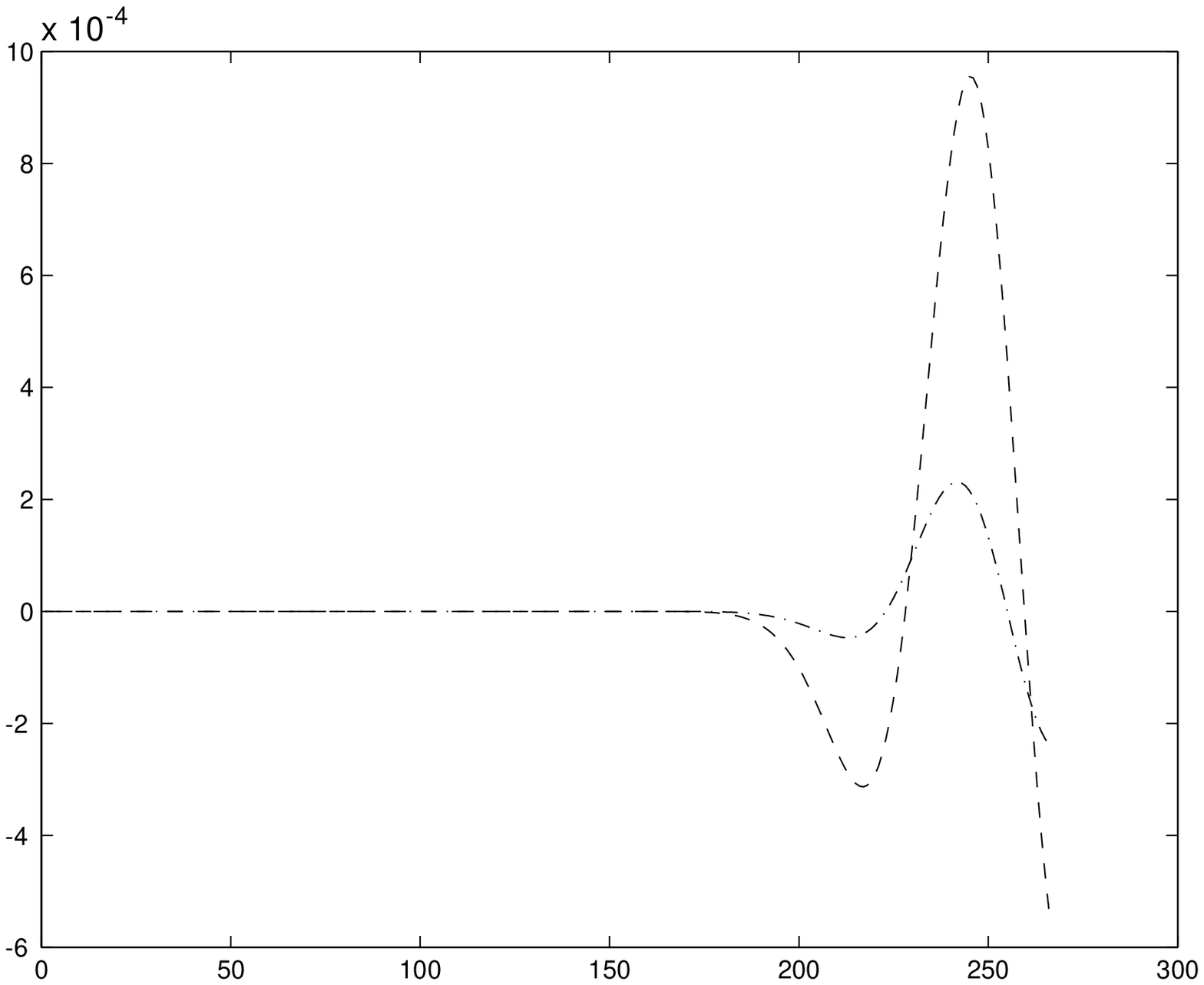}
  \caption{Pressure fields for the reference solution and the PML solution near the upperleft corner (left) and in the PML (right)  for an oblique flow $\ub=\vb=0.33$ vs. time steps}
  \label{fig:oblique}
\end{figure}

\begin{figure}
  \centering
  \includegraphics[width=0.48\textwidth,height=0.35\textwidth]
    {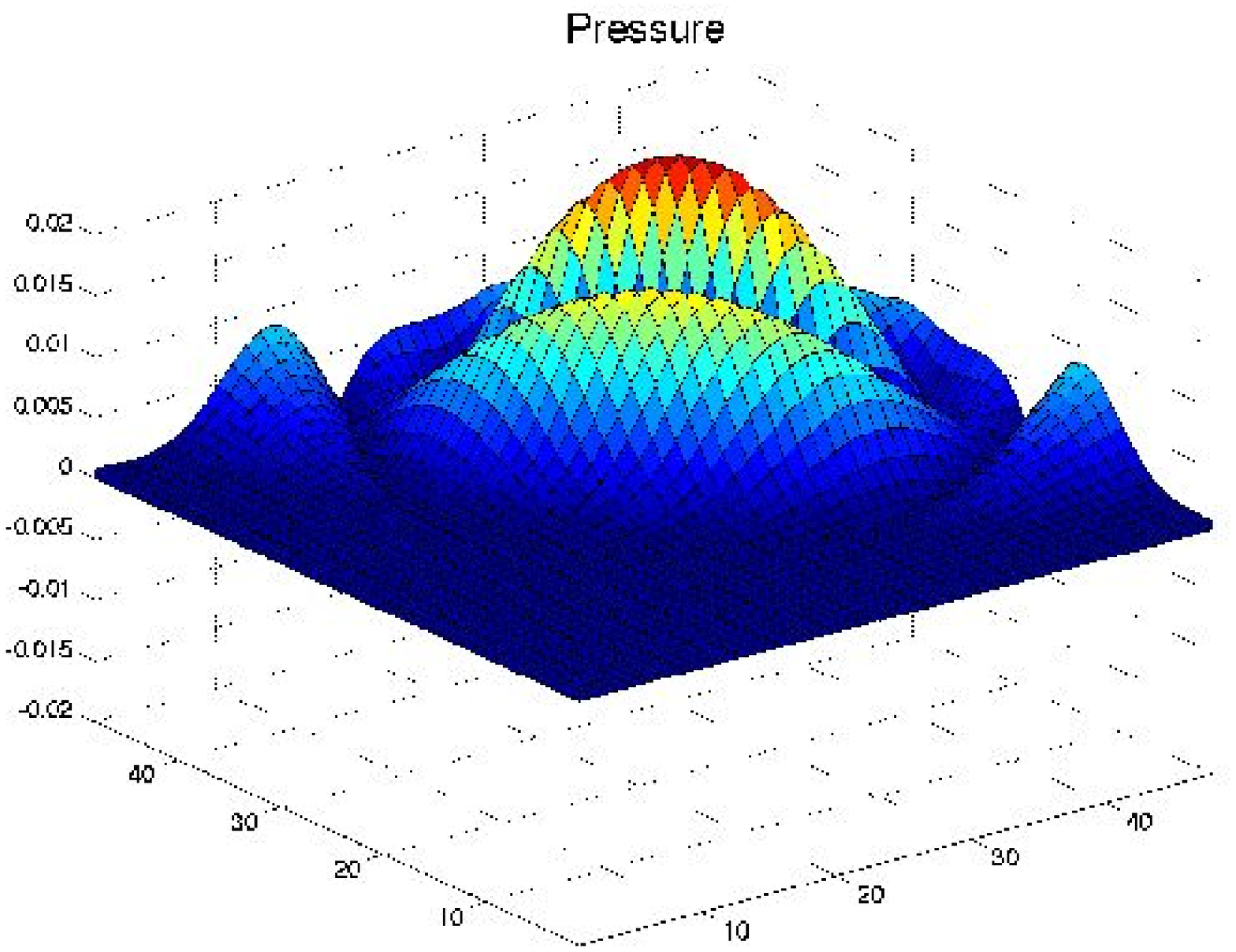}
  \includegraphics[width=0.48\textwidth,height=0.35\textwidth]
    {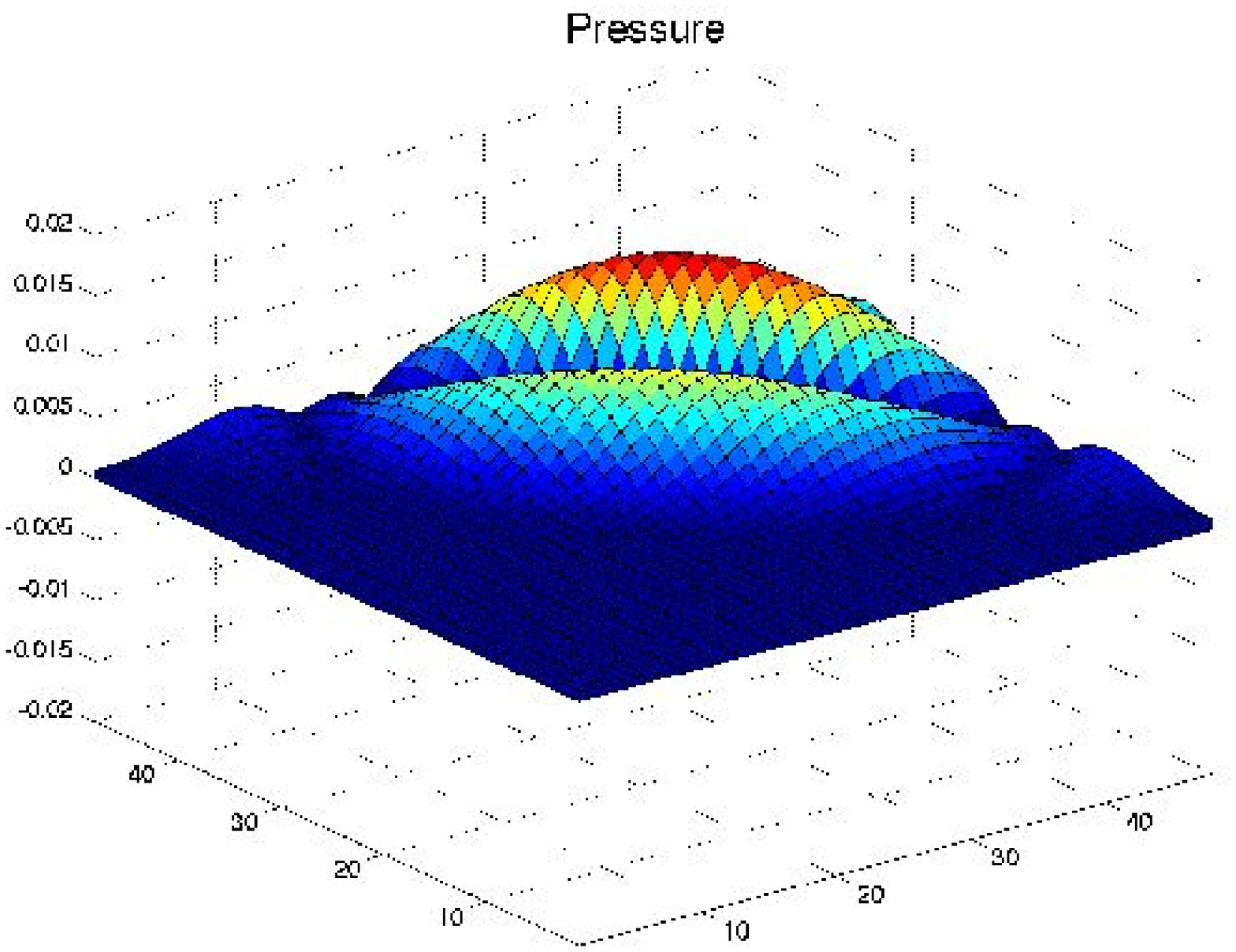}
  \includegraphics[width=0.48\textwidth,height=0.35\textwidth]
    {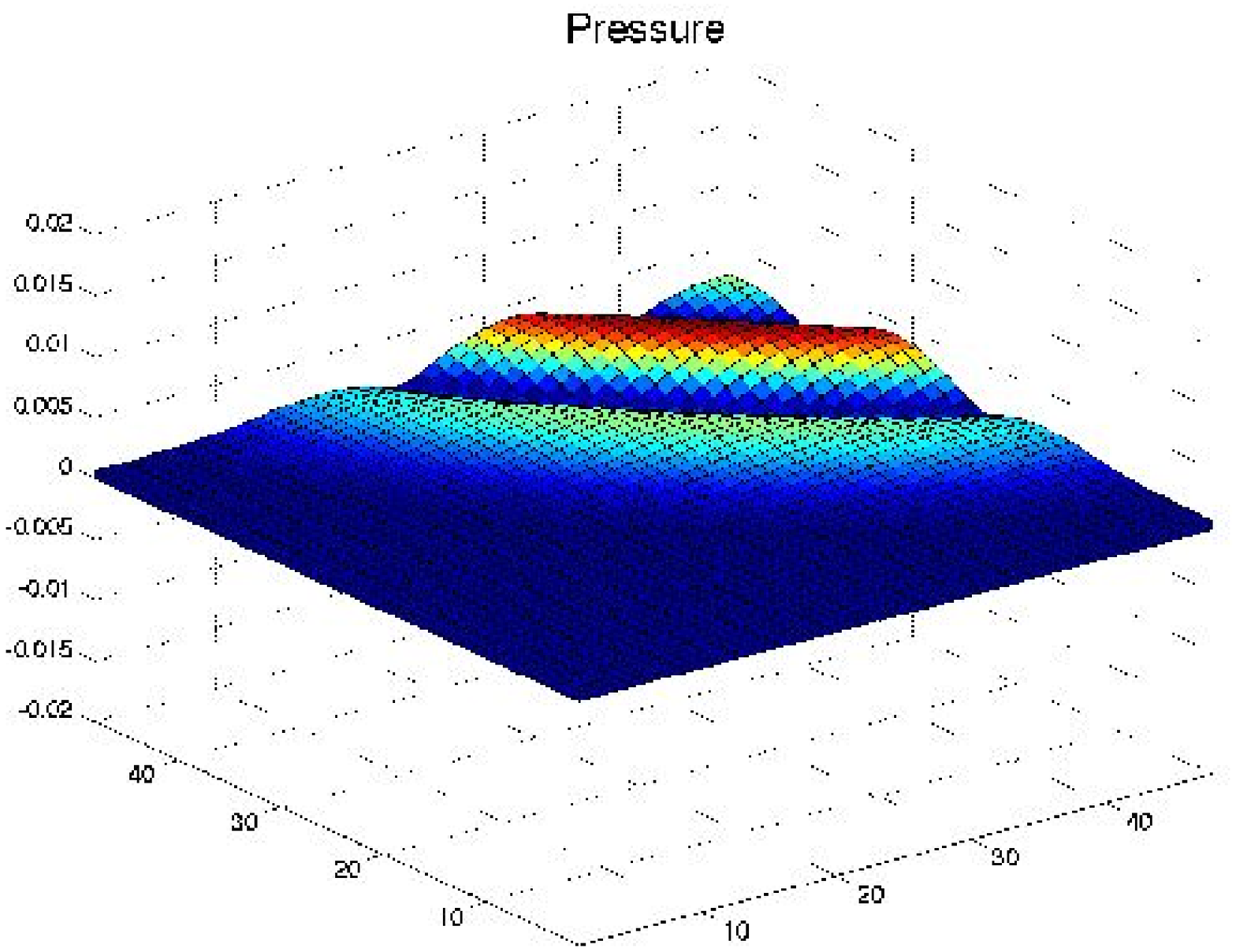}
  \includegraphics[width=0.48\textwidth,height=0.35\textwidth]
    {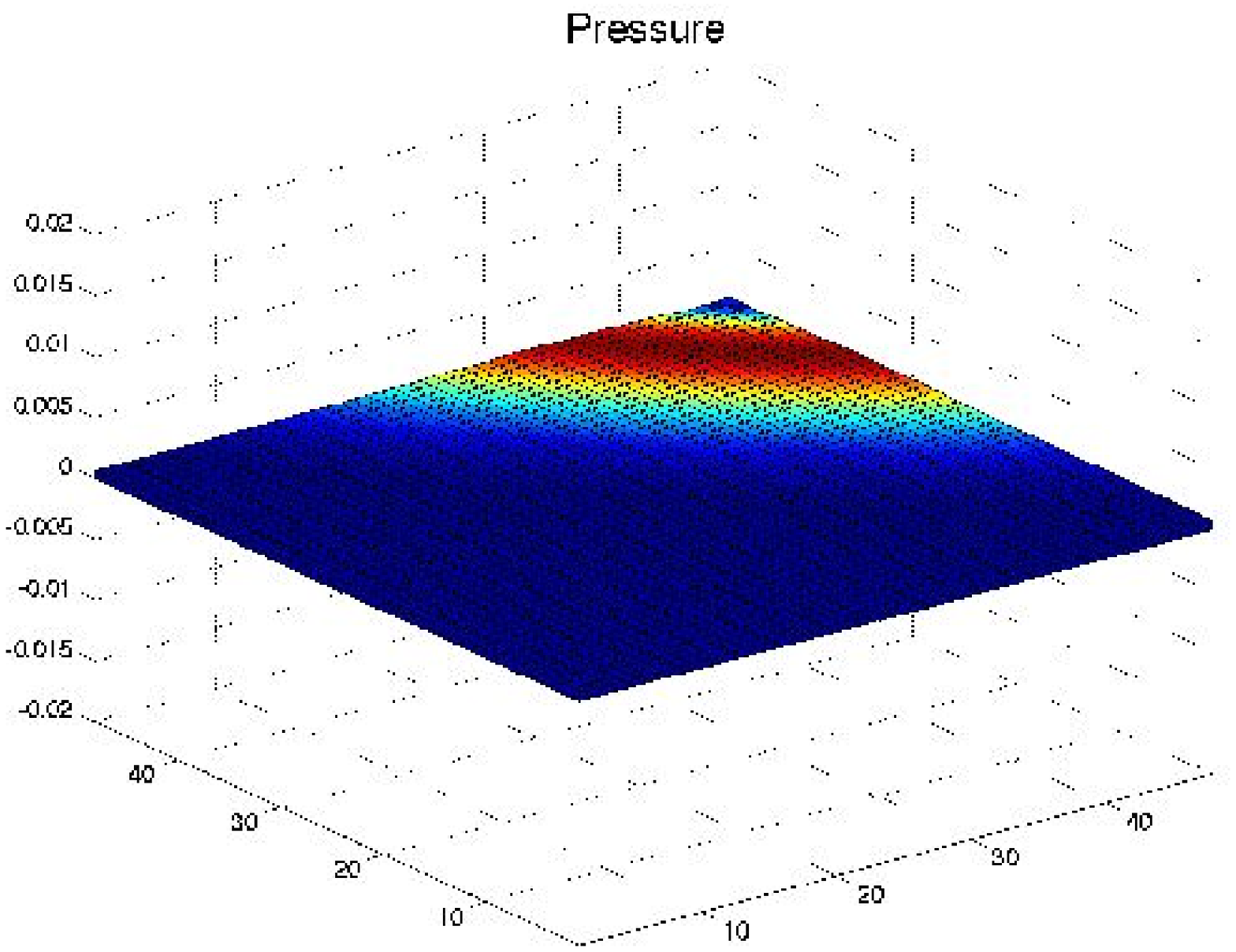}
  \caption{Pressure field for an oblique velocity $\ub=\vb=0.9$ at successive time steps}
  \label{fig:obliquecartoon}
\end{figure}

\begin{figure}
  \centering
  \includegraphics[width=0.96\textwidth,height=0.65\textwidth]
    {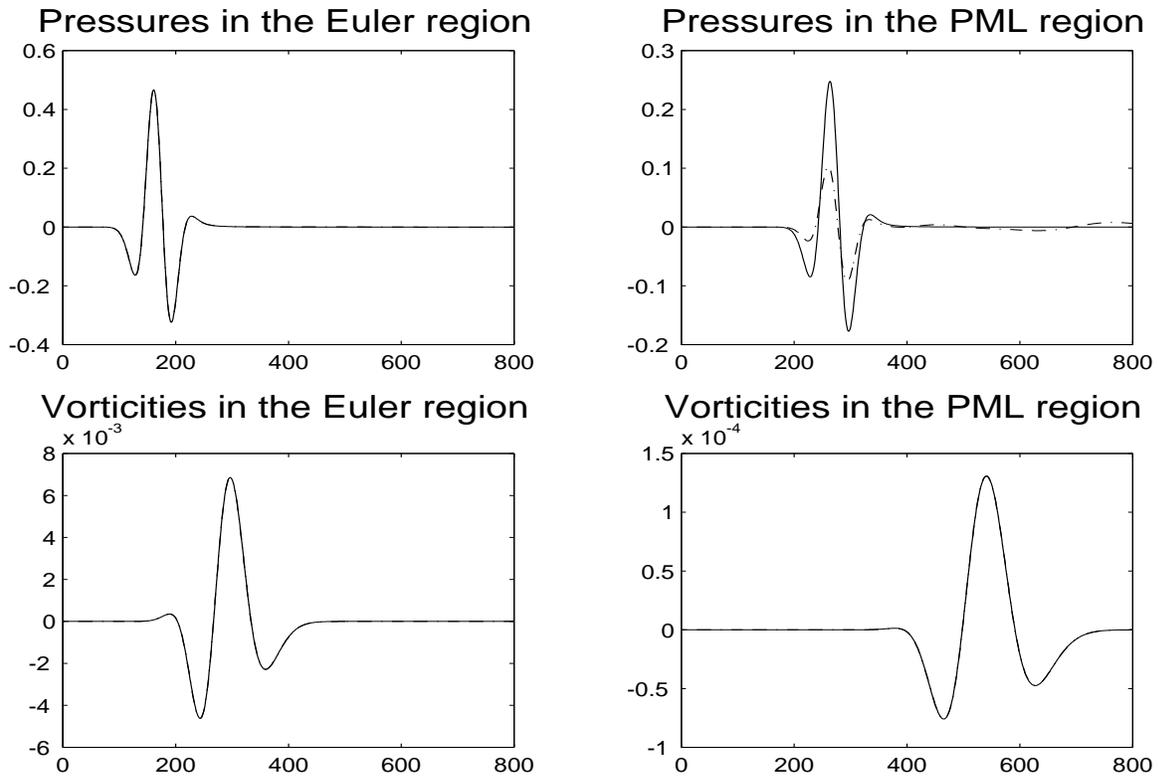}
  \caption{Reference and ``PML'' solutions in the Euler and PML regions vs. time steps ($\ub=200$, $\vb=100$)}
  \label{fig:omegauntouched}
\end{figure}

\begin{figure}
  \centering
  \includegraphics[width=0.48\textwidth,height=0.35\textwidth]
    {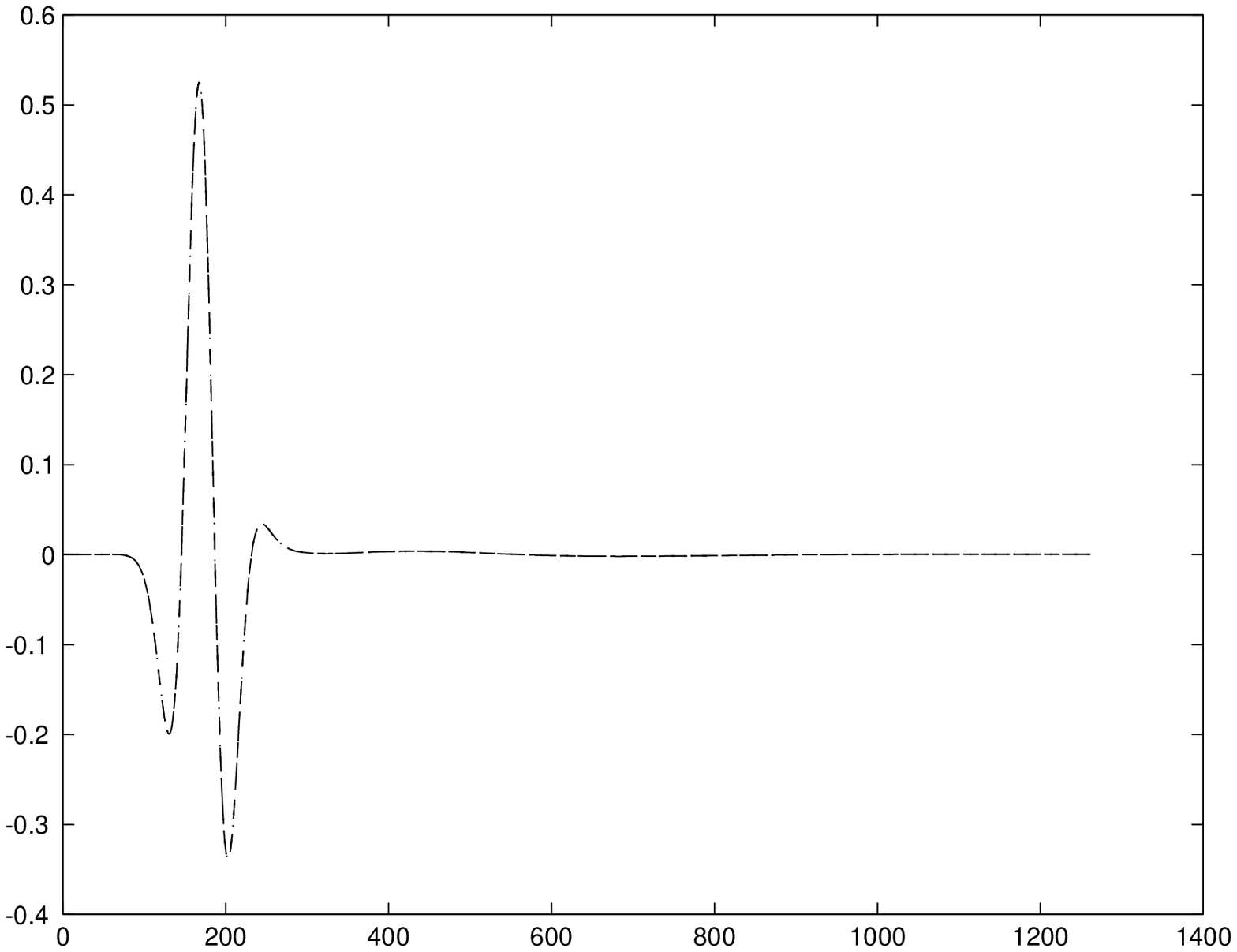}
  \includegraphics[width=0.48\textwidth,height=0.35\textwidth]
    {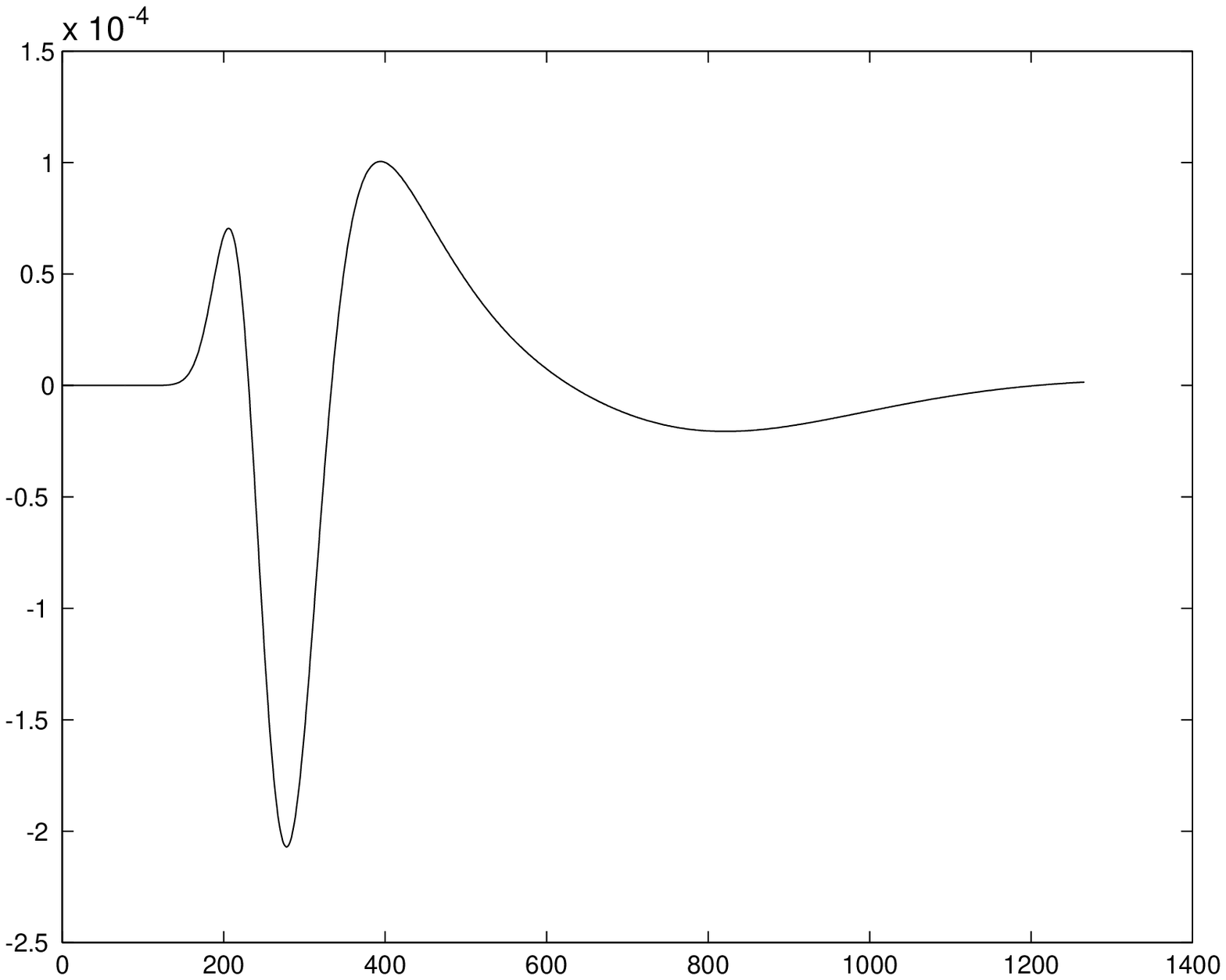}
  \caption{Pressure field (left) and error on the pressure (right) near the upperleft corner  for  an oblique velocity $\ub=\vb=0.9$ vs. time steps}
  \label{fig:oblique270}
\end{figure}

\begin{figure}
  \centering
  \includegraphics[width=0.48\textwidth,height=0.35\textwidth]
    {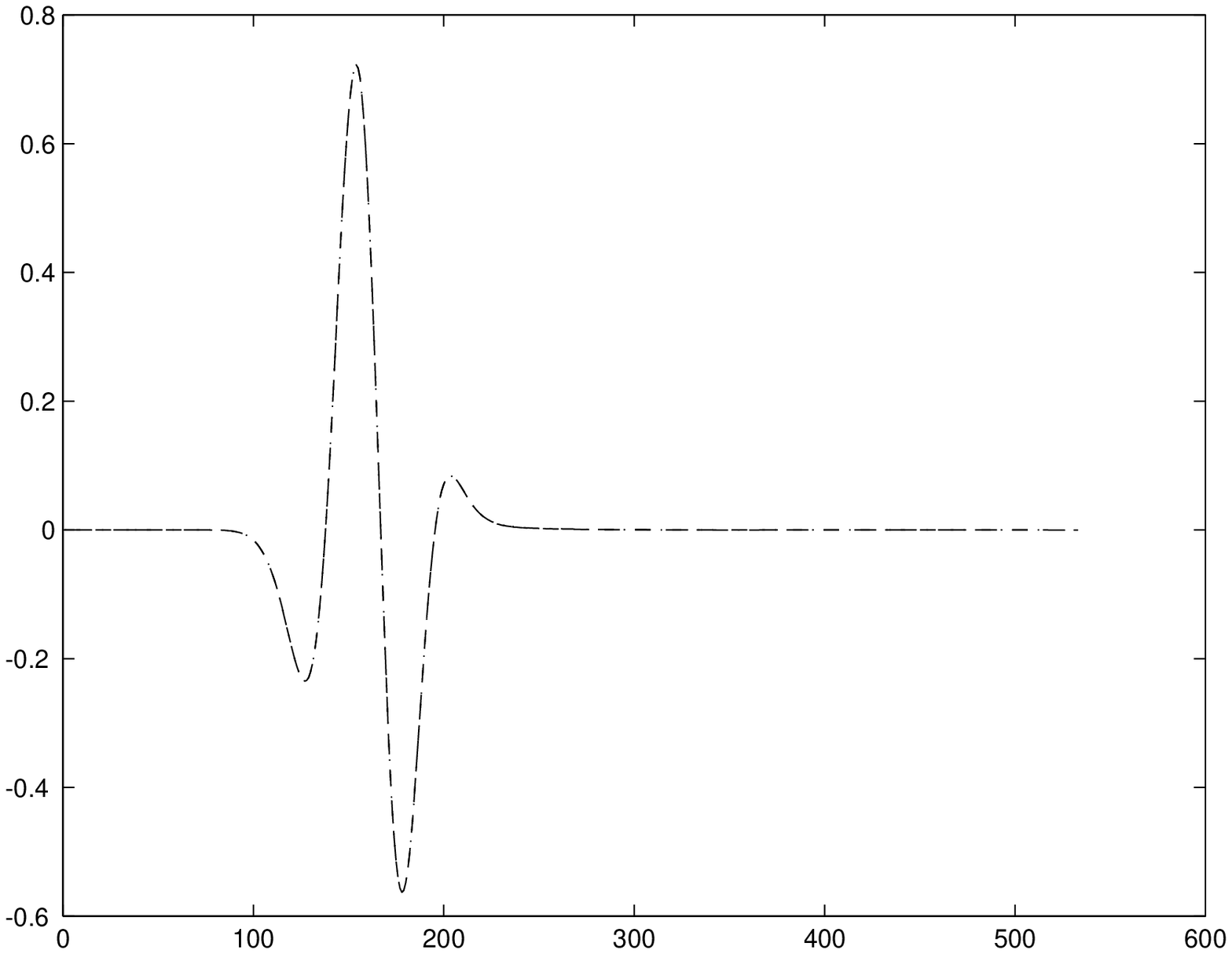}
  \includegraphics[width=0.48\textwidth,height=0.35\textwidth]
    {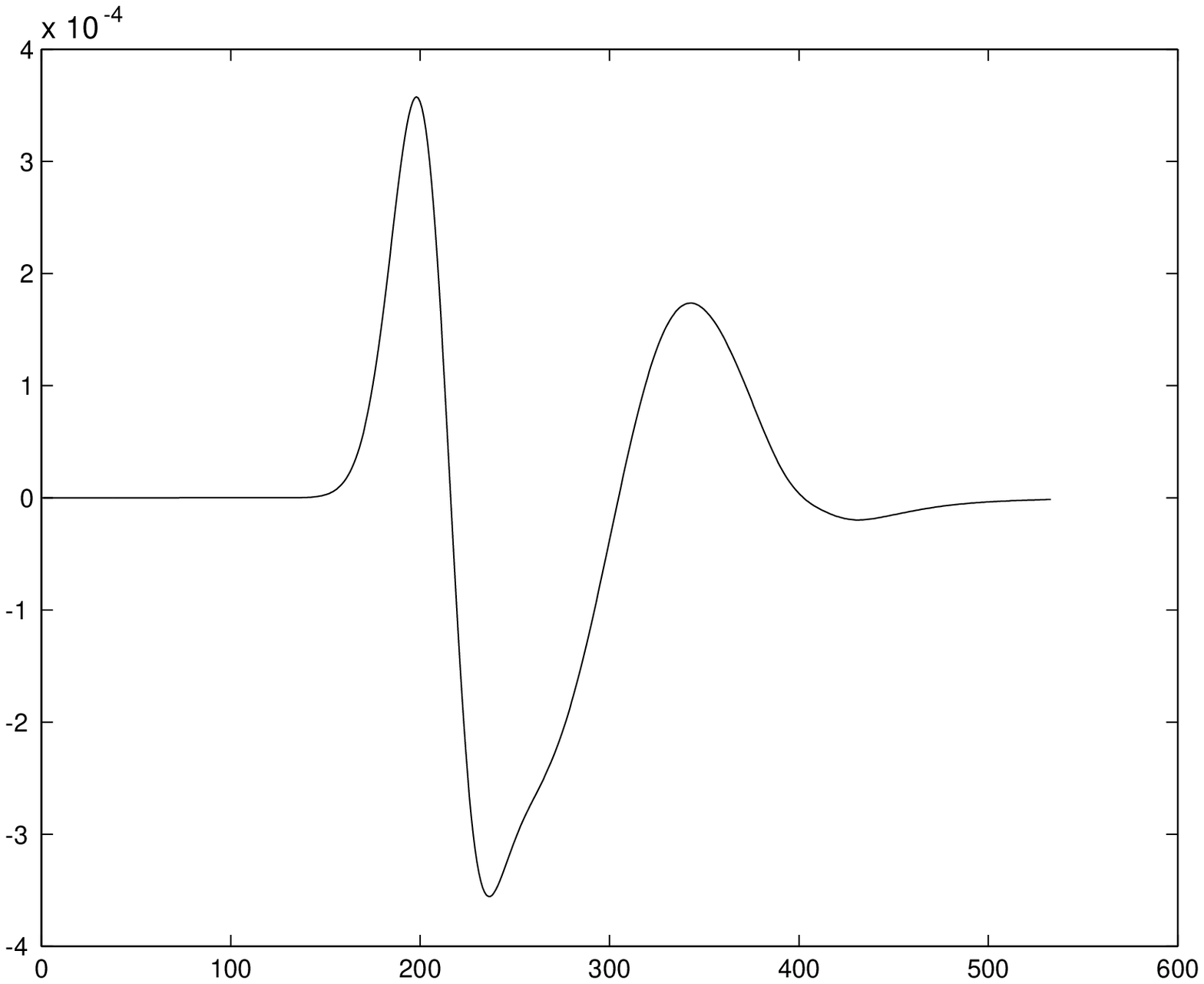}
  \caption{Pressure field (left) and error on the pressure (right) near the upperleft corner for a horizontal flow ($\ub=0.33,\,\vb=0$) in a duct vs.time steps}
  \label{fig:horizontal}
\end{figure}

We have taken $\cb=300$ and $\rhob=1$. The 2D linearized Euler equations are discretized on a uniform staggered grid using a Yee Scheme and a CFL equals to $0.3$. The convective derivatives are discretized using an upwind scheme both in the Euler region and in the PMLs. The computational domain is the square $[0,1.2]\times [0,1.2]$ and except in table~\ref{tab:errvswidthandsigma}, PMLs have a width $\delta=0.9$ corresponding to $n_\delta=38$ grid points. The damping parameter $\sigma$ depends on the coordinate normal to the interface (say $x$): $\sigma(x)=\sigma_{pml}\,x^2/\delta^3$ where $\sigma_{pml}$ is a positive constant. The initial solutions are zero. Let $f(t,x,y)=(1-2\pi^2(f_c t -1)^2) e^{-\pi^2(f_c t -1)^2}\delta_M(x,y)$ for $t<T_s$ and zero for $t>T_s$ with $T_s=0.05$, $f_c=4/T_s$ and $\delta_M$ is the Dirac mass located in the middle of the computational domain. The right handside was $f(t,x,y)$ on all three equations of system (\ref{eq:euler}) except for figure~\ref{fig:oblique} where it was zero on the velocity components. The PML solution is compared with a reference solution that is computed on a much larger domain. In figure~\ref{fig:oblique}, the pressure for both solutions are plotted as a function of time 4 points from the upperleft corner of the domain (left figure) and inside the PML (right figure). The velocity field is $\ub=\vb=1/3$. In the Euler region, both curves are nearly identical. In the PML, we see the damping of the PML solution. Of course, for the reference solution, this corresponds to an Euler region and there is no damping. In Figure~\ref{fig:omegauntouched} we have the same plot for the pressure and the vorticity as well, in this case $\ub=200$ and $\vb=100$. In agreement with the construction of the PML, we see that the vorticity is not damped at all in the PML region. The vorticity in the PML region equals that of the reference solution. Indeed, in the construction of the PML in \eqref{eq:eulerpressure} only the wave operator ${\mathcal L}$ is ``pml''-ized and the transport operator ${\mathcal G}$ is not modified. In figure~\ref{fig:obliquecartoon}, we show the pressure at different times of the computation for an oblique velocity $\ub=\vb=270$. For the same computation, pressure near the upperleft corner is shown on Figure~\ref{fig:oblique270}. Figure~\ref{fig:horizontal} is a similar figure for a horizontal flow in a duct.\\

In table~\ref{tab:errvswidthandsigma}, we study the influence of the parameters of the PML on the error between the reference solution and the pml solution. We see that with a generous PML ($38$ grid points), the error is indeed very small. With $18$ grid points, the error is small and with $8$ points results are not satisfactory even when using various parameters $\sigma_{pml}$. It is worth noticing that due to the use of an upwind scheme for the transport operator, the error on the vorticity is equal to the machine accuracy for any layer parameters. \\

\begin{table}[!h]
       \caption{Relative errors in percentage for different PML parameters, ($\ub=200$, $\vb=100$)}
       \label{tab:errvswidthandsigma}
     \begin{center}
     \begin{small}
       \begin{tabular}[c]{|c|r|r|r|r|r|r|}\hline
       &\multicolumn{6}{|c|}{Relative error in percentage vs. ($n_\delta$, $\sigma_{pml}$)} \\
\cline{2-7}
     variable & \multicolumn{1}{|c|}{($38$, $40$)}
              &\multicolumn{1}{|c|}{($18$, $40$)} 
              &\multicolumn{1}{|c|}{($18$, $80$)} 
              &\multicolumn{1}{|c|}{($8$, $40$)}
              &\multicolumn{1}{|c|}{($8$, $80$)}
              &\multicolumn{1}{|c|}{($8$, $160$)}
    \\ \hline
         $p$ & 0.12 & 1.0  & 8.0  & 10.0& 6.0  & 1.0 \\ \hline
         $u$ & 0.25 & 1.2 & 0.4 & 16.3 & 0.4 & 0.3 \\ \hline
         $v$ & 0.1  & 0.1 & 10.0  & 20.8 & 0.4 & 10.0 \\ \hline
    $\omega$ & $2.10^{-14}$ & $3.10^{-14}$ & $2.10^{-14}$  &$2.10^{-14}$  & $4.10^{-14}$  & $2.10^{-14}$    \\ \hline
       \end{tabular}
     \end{small}
    \end{center}
\end{table}

 The long time stability of the PML was assessed by computing on time intervals five times longer than those used for generating the figures. No instability was observed for various flows.

\section{Conclusion}\label{sec:conclusion}
The first PML model proposed in \S~\ref{sec:firstpml} is obtained by using the Smith factorization of the Euler equations and a PML for the advective wave equation. This method can be applied to other systems of partial differential equations (free-surface flow, anisotropic elasticity, $\ldots$). The second PML model we have proposed for the Euler linearized equations are based on the PML for the advective wave equation. Thus, the PML for Euler inherits the properties from the latter. This second model was implemented and numerical results illustrate the efficiency of the approach.

\bibliography{paperSysteme}

\end{document}